\documentclass[11pt]{article}
\usepackage{epsfig}
\usepackage{amssymb}
\usepackage{amscd}
\usepackage[all]{xy}
\usepackage{epsf}
\usepackage{amsfonts,amssymb,amsthm,amsmath}

\input cyracc.def

\long\def\comment#1\endcomment{\relax}

\newcounter{subsubsubsection}
\newcounter{subsubsubsubsection}

\makeatletter \@addtoreset{subsubsubsection}{subsubsection}
\@addtoreset{subsubsubsubsection}{subsubsubsection}

\makeatother

\DeclareMathSizes{11.1}{10}{8}{6}

\DeclareMathOperator{\Hom}{Hom}

\newtheorem*{theorem*}{Theorem}

\newtheorem*{conjecture}{Conjecture}

\newtheorem*{lemma}{Lemma}

\newtheorem*{proposition}{Proposition}

\theoremstyle{remark}

\newtheorem*{remark}{Remark}

\newtheorem*{example}{Example}

\theoremstyle{definition}

\newcommand{\hdot}{{\:\protect\raisebox{3pt}{\text{\protect\circle*{1.5}}}}}

\newcommand{\mb}{\hdot}

\newcommand{\Ext}{\mathrm{Ext}}

\newcommand{\g}{\mathfrak{g}}
\newcommand{\Sym}{\mathrm{Sym}}

\newcommand{\K}{\mathrm{K}}

\newcommand{\Der}{\mathrm{Der}}
\newcommand{\Inn}{\mathrm{Inn}}

\newcommand{\xxx}{\otimes\mathbb{C}[[\hbar]]}

\newcommand{\Hoch}{\mathrm{Hoch}}
\newcommand{\poly}{{poly}}

\title{ {\huge  Koszul duality in deformation quantization, I}}
\author{{\LARGE Boris Shoikhet}}
\date{}

\begin{document}\maketitle

\comment

{\font\tcyr=wncyi10

\tcyr\cyracc

\font\scyr=wncyr10

\scyr\cyracc

\hbox to\textwidth{\hfil\parbox{80mm}{\tcyr{Zachem ya algebry ne
znal podi po{\u\i}mi se{\u\i}chas}\\{\tcyr Tako{\u\i} polezny{\u\i}
delovo{\u\i} predmet...}}} \vspace{2mm}

\hbox to\textwidth{\hfil\parbox{30mm}{\scyr{M.Shcherbakov}}}}

\vspace{1cm}

\endcomment
\begin{abstract}
Let $\alpha$ be a polynomial Poisson bivector on a
finite-dimensional vector space $V$ over $\mathbb{C}$. Then
Kontsevich [K97] gives a formula for a quantization $f\star g$ of
the algebra $S(V)^*$. We give a construction of an algebra with the
PBW property defined from $\alpha$ by generators and relations.
Namely, we define an algebra as the quotient of the free tensor
algebra $T(V^*)$ by relations $x_i\otimes x_j-x_j\otimes
x_i=R_{ij}(\hbar)$ where $R_{ij}(\hbar)\in T(V^*)\otimes\hbar
\mathbb{C}[[\hbar]]$, $R_{ij}=\hbar
\Sym(\alpha_{ij})+\mathcal{O}(\hbar^2)$, with one relation for each
pair of $i,j=1...\dim V$. We prove that the constructed algebra
obeys the PBW property, and this is a generalization of the
Poincar\'{e}-Birkhoff-Witt theorem. In the case of a linear Poisson
structure we get the PBW theorem itself, and for a quadratic Poisson
structure we get an object closely related to a quantum $R$-matrix
on $V$. At the same time we get a free resolution of the deformed
algebra (for an arbitrary $\alpha$).

The construction of this PBW algebra is rather simple, as well as
the proof of the PBW property. The major efforts should be
undertaken to prove the conjecture that in this way we get an
algebra isomorphic to the Kontsevich star-algebra.
\end{abstract}

\section{The main construction}
\subsection{}
First of all, recall here the Stasheff's definition of the
Hochschild cohomological complex of an associative algebra $A$.

Consider the shifted vector space $W=A[-1]$, and the cofree
coassociative coalgebra $C(W)$ (co)generated by $W$. As a graded
vector space, $C(W)=T(A[1])$, the free tensor space. The coproduct
is:

\begin{equation}\label{eq1}
\begin{aligned}
\ &\Delta(a_1\otimes a_2\otimes\dots\otimes a_k)=\sum_{i=1}^{k-1}
(a_1\otimes\dots \otimes a_i)\bigotimes (a_{i+1}\otimes\dots\otimes
a_k)
\end{aligned}
\end{equation}

Consider the Lie algebra $CoDer(C(A[1]))$ of all coderivations of
this coalgebra. As the coalgebra is free, any coderivation $D$ (if
it is graded) is uniquely defined by a map $\Psi_D\colon A^{\otimes
k}\to A$, and the degree of this coderivation is $k-1$ (in
conditions that $A$ is not graded). The bracket
$[\Psi_{D_1},\Psi_{D_2}]$ is again a coderivation. Define the
Hochschild Lie algebra as $\Hoch^\mb(A)=CoDer^\mb(C(A[1]))$. To
define the complex structure on it, consider the particular
coderivation $D_m$ of degree +1 from the product $m\colon A^{\otimes
2}\to A$, which is the product in the associative algebra $A$. The
condition $[D_m,D_m]=0$ is equivalent to the associativity of $m$.
Define the differential on $CoDer^\mb(C(A[1]))$ as
$d(\Psi)=[D_m,\Psi]$. In this way we get a dg Lie algebra. The
differential is called the Hochschild differential, and the bracket
is called the Gerstenhaber bracket. The definition of these
structures given here is due to J.Stasheff.

\subsection{The explicit definition}
Here we relate the Stasheff's definition of the Hochschild
cohomological complex with the usual one.

The concept of a coderivation of a (co)free coalgebra is dual to the
concept of a derivation of a free algebra. Let $L$ be a vector
space, and let $T(L)$ be the free tensor algebra generated by the
vector space $L$. Let $D\colon T(L)\to T(L)$ be a derivation, then
it is uniquely defined by its value $D_L\colon L\to T(L)$ on the
generators, and any $D_L$ defines a derivation $D$ of the free
algebra $T(L)$. If we would like to consider only graded
derivations, we restrict ourselves by the maps $D_L\colon L\to
L^{\otimes k}$ for $k\ge 0$.

Dually, a coderivation $D$ of the cofree coalgebra $C(P)$
cogenerated by a vector space $P$ is uniquely defined by the
restriction to cogenerators, that is, by a map $D_P\colon C(P)\to
P$, or, if we consider the graded coderivations, the map $D_P$ is a
map $D_P\colon P^{\otimes k}\to P$ for $k\ge 0$.

In our case of the definition of the cohomological Hochschild
complex of an associative algebra $A$, we have $P=A[1]$. Then the
coderivations of the grading $k$ form the vector space
$\Hoch^k(A)=\Hom(A^{\otimes (k+1)},A)$, $k\ge -1$. Now we can deduce
the differential and the Gerstenhaber bracket from the Stasheff's
construction. The answer is the following:

For $\Psi\in\Hom(A^{\otimes k},A)$ the cochain $d\Psi\in
\Hom(A^{\otimes (k+1)},A)$ is given by the formula:

\begin{equation}\label{eq20}
\begin{aligned}
\ &d\Psi(a_0\otimes\dots\otimes a_k)=a_0\Psi(a_1\otimes\dots\otimes a_k)+\\
&+\sum_{i=0}^{k-1}(-1)^{i+1}\Psi(a_0\otimes\dots\otimes
a_{i-1}\otimes(a_ia_{i+1})\otimes a_{i+2}\otimes\dots\otimes a_k)+\\
&+(-1)^{k+1}\Psi(a_0\otimes\dots\otimes a_{k-1})a_k
\end{aligned}
\end{equation}

For $\Psi_1\in\Hom(A^{\otimes(k+1)},A)$ and
$\Psi_2\in\Hom(A^{\otimes (l+1)},A)$ the bracket
$[\Psi_1,\Psi_2]=\Psi_1\circ\Psi_2-(-1)^{kl}\Psi_2\circ\Psi_1$ where

\begin{equation}\label{eq21}
\begin{aligned}
\ & (\Psi_1\circ\Psi_2)(a_0\otimes\dots\otimes a_{k+l})=\\
&\sum_{i=0}^k(-1)^{il}\Psi_1(a_0\otimes\dots\otimes
a_{i-1}\otimes\Psi_2(a_i\otimes\dots\otimes a_{i+l})\otimes
a_{i+l+1}\otimes\dots\otimes a_{k+l})
\end{aligned}
\end{equation}

\subsection{The (co)bar-complex}
Here we recall the definition of the (co)bar-complex of an
associative (co)algebra. When the (co)algebra contains (co)unit, the
(co)bar-complex is acyclic, and when the (co)algebra is the kernel
of the augmentation, this concept is closely related to the Koszul
duality.

Let $A$ be an associative algebra. Then its bar-complex is
$$
\dots\rightarrow A^{\otimes 3}\rightarrow A^{\otimes 2}\rightarrow
A\rightarrow 0
$$
where $\deg A^{\otimes k}=-k+1$, and the differential $d\colon
A^{\otimes k}\to A^{\otimes (k-1)}$ is given as follows:

\begin{equation}\label{eq23}
d(a_1\otimes\dots\otimes a_k)=(a_1a_2)\otimes a_3\otimes\dots\otimes
a_k-a_1\otimes (a_2a_3)\otimes\dots\otimes a_k+\dots
+(-1)^{k}a_1\otimes\dots\otimes a_{k-2}\otimes (a_{k-1}a_k)
\end{equation}
If the algebra $A$ has unit, the bar-complex of $A$ is acyclic in
all degrees. Indeed, the map
$$
a_1\otimes \dots\otimes a_k\mapsto 1\otimes a_1\otimes \dots\otimes
a_k
$$
is a contracting homotopy.

Suppose now that the algebra $A$ does not contain unit, and $A=B^+$
is the kernel of an augmentation map $\varepsilon\colon
B\to\mathbb{C}$. (The map $\varepsilon$ is a surjective map of
algebras, in particular, it maps $1$ to $1$). Then the cohomology of
the bar-complex of $A$ is equal to the dual space
$\Ext^\mb_{B-Mod}(\mathbb{C},\mathbb{C})$.

Indeed, for any $B$-module $M$, we have the following free
resolution of $M$:

\begin{equation}\label{eq25}
\dots B\otimes \overline{B}\otimes\overline{B}\otimes M\rightarrow
B\otimes\overline{B}\otimes M\rightarrow B\otimes M\rightarrow
M\rightarrow 0
\end{equation}
with the differential analogous to the bar-differential.

Consider the case $M=\mathbb{C}$. We can compute
$\Ext_{B-Mod}^\mb(\mathbb{C},\mathbb{C})$ using this resolution. In
the answer we get the cohomology of the complex dual to the
bar-complex of $A=B^+$.

The complex dual to the bar-complex of $A$ is the {\it
cobar-complex} for the coalgebra $A^*$. This cobar-complex is an
associative dg algebra, and it is a free algebra, which by previous
is a free resolution of the algebra
$\Ext_{B-Mod}(\mathbb{C},\mathbb{C})$. For the sequel we write down
explicitly the cobar-complex of a coassociative coalgebra $Q$:

\begin{equation}\label{eq26}
0\rightarrow Q\rightarrow Q\otimes Q\rightarrow Q\otimes Q\otimes
Q\rightarrow\dots
\end{equation}
and the differential $\delta Q^{\otimes k}\to Q^{\otimes (k+1)}$ is
\begin{equation}\label{eq27}
\delta(q_1\otimes\dots\otimes q_k)=(\Delta q_1)\otimes
q_2\otimes\dots\otimes q_k-q_1\otimes(\Delta q_2)\otimes
\dots\otimes q_k+\dots+(-1)^{k-1}q_1\otimes \dots\otimes
q_{k-1}\otimes (\Delta q_k)
\end{equation}
where $\Delta\colon Q\to Q^{\otimes 2}$ is the coproduct.

In the case when $B=S(V)$ is the symmetric algebra,
$\Ext_{B-Mod}(\mathbb{C},\mathbb{C})$ is the exterior algebra
$\Lambda (V^*)=S(V^*[-1])$, and vise versa. In this way, we get a
free resolution of the symmetric (exterior) algebra.

\begin{example}
Here we construct explicitly the free cobar-resolution
$\mathcal{R}^\mb$ of the algebra $\mathbb{C}[x_1,x_2]$ of
polynomials on two variables. As a graded algebra, $\mathcal{R}^\mb$
is the free algebra
$\mathcal{R}^\mb=\mathrm{Free}(x_1,x_2,\xi_{12})$ where $\deg
x_1=\deg x_2=0$, $\deg\xi_{12}=-1$. The differential is $0$ on
$x_1,x_2$, $d(\xi_{12})=x_1\otimes x_2-x_2\otimes x_1$, and
satisfies the graded Leibniz rule. In degree 0 we have the tensor
algebra $T(x_1,x_2)$, differential is 0 on degree 0 (there are no
elements in degree 1). In degree -1, a general element is a
non-commutative word in $x_1,x_2,\xi_{12}$ in which $\xi_{12}$
occurs exactly one time. For example, it could be a word $x_2\otimes
x_1\otimes x_2\otimes \xi_{12}\otimes x_1\otimes x_2$. The image of
the differential is then exactly the two-sided ideal in the tensor
algebra $T(x_1,x_2)$ generated by $x_1\otimes x_2-x_2\otimes x_1$.
Then, the 0-th cohomology is $\mathbb{C}[x_1,x_2]$. It follows from
the discussion above that all higher cohomology is 0.
\end{example}

\subsection{The main construction} Here we construct a quasi-isomorphic map of dg Lie
algebras $\Phi\colon \Hoch^\mb(S(V))/\mathbb{C}\to
\Der(CoBar^\mb(S(V^*)^+))/\Inn (CoBar^\mb(S(V^*)^+))$.

Let $\Psi\in\Hom ((S(V))^{\otimes k},S(V))$ be a $k$-cochain. Denote
$V^*=W$, then we can consider $\Psi$ the corresponding cochain in
$\Hom(S(W),(S(W))^{\otimes k})$. Here we consider $S(W)$ as {\it
coalgebra}. Then this cochain may be considered as a derivation in
$\Der(CoBar^\mb(S(W)))$. We would like to attach to it a derivation
in $\Der(CoBar^\mb(S(W)^+))$, maybe modulo an inner derivation. So,
we would like to show that there exist a map $\Phi\colon
\Der(CoBar^\mb(S(W)))\to\Der(CoBar^\mb(S(W)^+))$ such that the
diagram

\begin{equation}\label{diagram1}
\xymatrix{ \Der(CoBar^\mb(S(W)))\ar[r]^{\delta}  \ar[d]^{\Phi}&
\Der(CoBar^\mb(S(W)))\ar[d]^{\Phi}\\
\Der(CoBar^\mb(S(W)^+))\ar[r]^{\delta} &\Der(CoBar^\mb(S(W)^+))}
\end{equation}
is commutative {\it modulo inner derivations} (here $\delta$ is the
cobar differential).

In the coalgebra $S(W)$ the coproduct is given by the formula
\begin{equation}\label{eq27}
\Delta(x_1\dots x_k)=1\otimes (x_1\dots x_k)+\sum_{\{i_1\dots
i_a\}\sqcup\{j_1\dots j_b\}=\{1\dots k\}}(x_{i_1}\dots
x_{i_a})\otimes(x_{j_1}\dots x_{j_b})+(x_1\dots x_k)\otimes 1
\end{equation}
and in the coalgebra $S(W)^+$ the coproduct is given by the same
formula without the first and the last summands, which contain 1's.

Therefore the projection $p\colon S(W)\to S(W)^+$ is a map of
coalgebras (dual to the imbedding of algebras), and the imbedding
$i\colon S(W)^+\to S(W)$ is {\it not}.

If $\Psi\colon S(W)\to S(W)^{\otimes k}$ is as above, we define
$(\Phi(\Psi))(\sigma)=p^{\otimes k}(\Psi(i(\sigma)))\in\Hom (S(W)^+,
(S(W)^+)^{\otimes k})$.

Now we check the commutativity of the diagram (\ref{diagram1})
modulo inner derivations. It is clear that

\begin{equation}\label{eq28}
(\Phi\circ\delta)(\sigma)-(\delta\circ\Phi)(\sigma)=p^{\otimes
k}(\Psi(1))\otimes\sigma\pm \sigma\otimes p^{\otimes k}(\Psi(1))
\end{equation}
which is an inner derivation $ad(p^{\otimes k}(\Psi(1)))$.

We have defined a map $\Phi_1\colon \Der(CoBar^\mb(S(W)))\to
\Der(CoBar^\mb(S(W)^+))/\Inn(CoBar^\mb(S(W)^+))$. The first dg Lie
algebra is clearly isomorphic to the Hochschild cohomological
complex of the algebra $S(V)$ modulo constants, and we can consider
the map $\Phi_1$ as a map
$$
\Phi_1\colon\Hoch^\mb(S(V))/\mathbb{C}\to\Der(CoBar^\mb(S(V^*)^+))/\Inn(CoBar^\mb(S(V^*)^+))
$$
Let us note that the cobar complex $CoBar^\mb(S(V^*)^+)$ is a free
resolution of the Koszul dual algebra $\Lambda(V^*)$.

Now we have the following result:

\begin{proposition}
The map $\Phi_1$ is a quasi-isomorphism of dg Lie algebras.
\end{proposition}
It is clear that $\Phi_1$ is a map of dg Lie algebras, one only
needs to proof that it is a quasi-isomorphism of complexes. Although
it is useful to have this Proposition in mind when reading Section
2, we will not use it. The proof will appear somewhere.

\section{Applications to deformation quantization}
\subsection{A lemma}
We start with the following lemma, which is a formal version of the
semi-continuity of cohomology of a complex depending on a parameter,
which says that in a "singular value" of the parameter the
cohomology may only raise:
\begin{lemma}
Let $\mathcal{R}^\mb$ be a $\mathbb{Z}_{\le 0}$-graded complex with
differential $d$, such that $H^i(\mathcal{R}^\mb)$ vanishes for all
$i\ne 0$. Consider
$\mathcal{R}^\mb_{\hbar}=\mathcal{R}^\mb\otimes\mathbb{C}[[\hbar]]$.
Let $d_h\colon \mathcal{R}^\mb_\hbar\to
\hbar\mathcal{R}^{\mb+1}_\hbar$ be a linear map of degree +1 such
that
$$
(d+d_\hbar)^2=0
$$
Then the cohomology $H^i_\hbar$ of the complex
$\mathcal{R}^\mb_\hbar$ with the differential $d+d_\hbar$ vanishes
for $i\ne 0$, and as a vector space,
$H^0_\hbar(\mathcal{R}^\mb_\hbar,d+d_\hbar)\simeq
H^0(\mathcal{R}^\mb,d)\otimes \mathbb{C}[[\hbar]]$.
\begin{proof}
Consider the filtration
$$
\mathcal{R}^\mb_\hbar\supset\hbar\mathcal{R}^\mb_\hbar\supset\hbar^2\mathcal{R}_\hbar\supset\dots
$$
of the complex $\mathcal{R}_\hbar^\mb$ with the differential
$d+d_\hbar$. Compute the cohomology of $\mathcal{R}_\hbar^\mb$ by
the spectral sequence corresponding to this filtration. The term
$E_0^{p,q}=\hbar^p\mathcal{R}^{p+q}_\hbar/\hbar^{p+1}\mathcal{R}^{p+q}_\hbar$,
and $d_\hbar$ acts by 0 on $E^{p,q}_0$. Therefore, the cohomology in
this term is the cohomology of the differential $d$ and is
$E_1^{p,-p}=\hbar^p H^0(\mathcal{R}^\mb,d)/\hbar^{p+1}
H^0(\mathcal{R}^\mb,d)$ and $E_1^{p,q}=0$ for $q\ne -p$. All higher
differentials are 0 by the dimensional reasons, and the spectral
sequence collapses in the term $E_1$. The spectral sequence clearly
converges to the cohomology of $\mathcal{R}_\hbar^\mb$.

Lemma is proven.
\end{proof}
\end{lemma}
\subsection{A proof of the classical Poincar\'{e}-Birkhoff-Witt theorem}
Let $\g$ be a Lie algebra. Its universal enveloping algebra
$\mathcal{U}(\g)$ is defined as the quotient-algebra of the tensor
algebra $T(\g)$ by the two-sided ideal generated by elements
$a\otimes b-b\otimes a-[a,b]$ for any $a,b\in \g$. The
Poincar\'{e}-Birkhoff-Witt theorem says that $\mathcal{U}(\g)$ is
isomorphic to $S(\g)$ as a $\g$-module. We suggest here a (probably
new) proof of this classical theorem, which certainly is not the
simplest one, but sheds some light on the cohomological nature of
the theorem.

Before starting with the proof, let us make some remark. Let us
generalize the universal enveloping algebra as follows. Consider the
tensor algebra $T(x_1,\dots,x_n)$ and its quotient $A_{c_{ij}^k}$ by
the two-sided ideal generated by the relations $x_i\otimes
x_j-x_j\otimes x_i-\sum_{k}c_{ij}^k x_k$, $1\le i<j\le n$, where
$c_{ij}^k$ are not supposed to satisfy the Jacobi identity
\begin{equation}\label{eq2.2.1}
\sum_a(c_{ij}^a c_{ak}^b+c_{jk}^a c_{ai}^b+c_{ki}^a c_{aj}^b)=0
\end{equation}
Then, if (\ref{eq2.2.1}) is not satisfied, the algebra
$A_{c_{ij}^k}$ is smaller than $S(x_1,\dots,x_k)$, that is, the
two-sided ideal, generated by the relations, is bigger than in the
Lie algebra case when (\ref{eq2.2.1}) is satisfied.

Now we pass to the proof. Let $\g$ be a Lie algebra. By the
discussion in Section 1.3, $CoBar^\mb(\Lambda^+(\g))$ is a free
resolution of the symmetric algebra $S(\g)$. Denote the
cobar-differential by $d$. Introduce in
$CoBar_\hbar^\mb=CoBar^\mb(\Lambda^+(\g))\otimes\mathbb{C}[[\hbar]]$
a new differential $d+d_\hbar$, where $d_\hbar\colon
CoBar_\hbar^\mb\to\hbar \cdot CoBar^{\mb+1}_\hbar$ comes from the
chain differential in the Lie homology complex
$\partial\colon\Lambda^i(\g)\to\Lambda^{i-1}(\g)$. We denote
$$
d_\hbar=\hbar\partial
$$
The equation $(d+d_\hbar)^2=0$ follows from the fact that the chain
Lie algebra complex is a dg coalgebra, and, therefore, its
cobar-complex is well-defined.

Now, by Lemma 2.1, the complex $CoBar^\mb_\hbar$ has only 0 degree
cohomology, which is isomorphic to
$H^0(CoBar^\mb(\Lambda^+(\g)))\otimes\mathbb{C}[[\hbar]]=S(\g)\otimes\mathbb{C}[[\hbar]]$
as a (filtered) vector space. On the other hand, we can compute 0-th
cohomology of $(CoBar_\hbar^\mb(\Lambda^+(\g)),d+d_\hbar)$ directly.
It is the quotient of the tensor algebra
$T(\g)\otimes\mathbb{C}[[\hbar]]$ by the two-sided ideal generated
by the relations $a\otimes b-b\otimes a-\hbar [a,b]$, $a,b\in\g$.

The specialization of the last isomorphism for $\hbar=1$ gives the
Poincar\'{e}-Birkhoff-Witt theorem.

\subsection{}
Consider the following sequence of maps:
\begin{equation}\label{eq2.3.1}
\begin{aligned}
\
&T_\poly(V^*)\xrightarrow{\mathcal{U}_S}\Hoch(S(V))\simeq\Der(CoBar(S(V^*)))\xrightarrow{\Phi_1}\\
&\xrightarrow{\Phi_1}\Der(CoBar(S^+(V^*)))/\Inn(CoBar(S^+(V^*)))
\end{aligned}
\end{equation}
Here the first map is the Kontsevich formality $L_\infty$ morphism
for the algebra $S(V)$, the second isomorphism follows from the
Stasheff's construction, and the third map is the map $\Phi_1$
defined in Section 1.4.

Apply now the composition (\ref{eq2.3.1}) to the vector space $V[1]$
instead of $V^*$.
\begin{lemma}
Let $V$ be a finite-dimensional vector space. Then there is a
canonical isomorphism of the graded Lie algebras $T_\poly(V^*)\simeq
T_\poly(V[1])$.
\begin{proof}
It is straightforward. The map maps $k$-polyvector field with
constant coefficients on $V^*$ to a $k$-linear function on $V[1]$,
and so on.
\end{proof}
\end{lemma}
\begin{remark}
The algebras $S(V)$ and $\Lambda(V^*[-1])$ are Koszul dual, and they
have isomorphic Hochschild comology with all structures (see [Kel]).
\end{remark}
Denote by $\K$ the correspondence $\K\colon T_\poly(V^*)\to
T_\poly(V[1])$ from Lemma. Let $\alpha$ be a polynomial Poisson
bivector on the space $V^*$. By the correspondence from Lemma, we
get a polyvector field $\K(\alpha)$ which in general is not a
bivector, but still satisfies the Maurer-Cartan equation
\begin{equation}\label{eq2.3.3}
[\K(\alpha),\K(\alpha)]=0 \end{equation}

Let us rewrite (\ref{eq2.3.1}) for $V[1]$:
\begin{equation}\label{eq2.3.4}
\begin{aligned}
\
&T_\poly(V[1])\xrightarrow{\mathcal{U}_\Lambda}\Hoch(\Lambda(V^*))\simeq\Der(CoBar(\Lambda(V)))\xrightarrow{\Phi_1}\\
&\xrightarrow{\Phi_1}\Der(CoBar(\Lambda^-(V)))/\Inn(CoBar(\Lambda^-(V)))
\end{aligned}
\end{equation}
 The composition
(\ref{eq2.3.4}) maps the polyvector $\hbar\K(\alpha)$ to a
derivation $d_\hbar$ of degree +1 in
$\Der(CoBar(\Lambda^-(V)))\xxx$, which satisfies the Maurer-Cartan
equation
\begin{equation}\label{eq2.3.2}
(d+d_\hbar)^2=0
\end{equation}
in $\Der/\Inn$, where $d$ is the cobar-differential.

Actually, (\ref{eq2.3.2}) is satisfied in
$\Der(CoBar(\Lambda^-(V)))\xxx$, not only in
$\Der(CoBar(\Lambda^-(V)))\xxx/\Inn(CoBar(\Lambda^-(V)))\xxx$.
Indeed, we suppose that $V$ is placed in degree 0, then
$CoBar(\Lambda^-(V))$ is $\mathbb{Z}_{\le 0}$-graded. Therefore, any
inner derivation has degree $\le 0$, while $d_\hbar$ has degree +1.
We have the following
\subsection{}
\begin{lemma}
Let $\alpha$ be a Poisson bivector on $V^*$, and let $\K(\alpha)$ be
the corresponding Maurer-Cartan polyvector of degree 1 in
$T_\poly(V[1])$. Then (\ref{eq2.3.4}) defines an $\hbar$-linear
derivation $d_\hbar$ of
$CoBar(\Lambda^-(V))\otimes\mathbb{C}[[\hbar]]$ of degree +1
corresponding to $\hbar\K(\alpha)$, such that
$$
(d+d_\hbar)^2=0
$$
where $d$ is the cobar-differential. Moreover, $d_\hbar$ obeys
\begin{equation}\label{eq2.3.5}
d_\hbar(\xi_i\wedge\xi_j)=\hbar\Sym(\alpha_{ij})+\mathcal{O}(\hbar^2)
\end{equation}
where $\xi_i\wedge\xi_j\in \Lambda^{2}(V)$, $\Sym(\alpha_{ij})\in
T(V)$ is the symmetrization, and $\{\xi_i\}$ is the basis in $V[1]$
dual to the basis $\{x_i\}$ in $V^*$ in which
$\alpha=\sum_{ij}\alpha_{ij}\partial_i\wedge\partial_j$
\begin{proof}
We only need to prove (\ref{eq2.3.5}), all other statements are
already proven. We prove it in details in Section 2.7.
\end{proof}
\end{lemma}
\subsection{}
Let $A$ be an $\hbar$-linear associative algebra which is the
quotient of the tensor algebra $T(V)\xxx$ of a vector space $V$ by
the two-sided ideal generated by relations
$$
x\otimes y-y\otimes x=R(x,y)
$$
for any $x,y\in V$, where $R(x,y)\in \hbar T(V)\xxx$. Consider the
following filtration:
\begin{equation}\label{eq2.5.1}
A\supset\hbar A\supset\hbar^2 A\supset\hbar^3 A\supset\dots
\end{equation}
This is clearly an algebra filtration: $(\hbar^k A)\cdot(\hbar^\ell
A)\subset\hbar^{k+\ell}A$. Consider the associated graded algebra
$\mathrm{gr} A$. We say that the algebra $A$ is a
Poincar\'{e}-Birkhoff-Witt (PBW) algebra if $\mathrm{gr}A\simeq
S(V)\xxx$ as a graded $\mathbb{C}[[\hbar]]$-linear algebra.

In general, $\mathrm{gr}A$ is less than
$S(V)\otimes\mathbb{C}[[\hbar]]$, it is a quotient of
$S(V)\otimes\mathbb{C}[[\hbar]]$. One can say that the PBW property
is equivalent to the property that the quotient-algebra has "the
maximal possible size".

\subsection{}
Let $\alpha$ be a polynomial Poisson bivector in $V^*$. In Sections
2.3 and 2.4 we constructed an $\hbar$-linear derivation $d_\hbar$ on
$CoBar(\Lambda^-(V))\xxx$ such that $(d+d_\hbar)^2=0$ where $d$ is
the cobar-differential. By Section 1.3, the cobar-complex
$CoBar(\Lambda^-(V))$ is a free resolution of the algebra $S(V)$, in
particular, the cohomology of $d$ does not vanish only in degree 0
where it is equal to $S(V)$. We are in the situation of Lemma 2.1.
In particular, the dg algebra $(CoBar(\Lambda^-(V))\xxx,d+d_\hbar)$
has only non-vanishing cohomology in degree 0, and this 0-degree
cohomology is an algebra, which is a PBW algebra by Lemma 2.1.
\begin{theorem*}
The construction above constructs from a Poisson polynomial bivector
$\alpha$ on $V^*$ an algebra $A_\alpha$ with generators
$x_1,\dots,x_n$ and relations $[x_i,x_j]=d_\hbar(\xi_i\wedge\xi_j)$.
This algebra is a PBW algebra. \qed
\end{theorem*}

\begin{conjecture}
The algebra $A_\alpha$ is isomorphic to the Kontsevich star-algebra
on $S(V)\xxx$ constructed from the Poisson bivector $\alpha$. (We
suppose that in the formality morphisms $\mathcal{U}_\Lambda\colon
T_\poly(V[1])\to\Hoch(\Lambda(V^*))$ in (\ref{eq2.3.4}), and
$\mathcal{U}_S\colon T_\poly(V^*)\to\Hoch(S(V))$ which is used in
the construction of the star-product, one uses the same propagator
in the definition of the Kontsevich integrals, see [K97]).
\end{conjecture}

In our approach, we lift this Conjecture on the level of complexes,
and get a diagram commutative in the homotopical category
$\mathcal{H}om_{dg}$ from [Sh]. The commutativity of this diagram
implies the Conjecture.

\subsection{An explicit formula}
One can write down explicitly the relations in the algebra
$A_\alpha$, in the terms of the Kontsevich integrals [K97]. For this
we need to find explicit formula for the $\hbar$-linear derivation
$d_\hbar$ in $CoBar(\Lambda^-(V))\xxx$. Here we suppose some
familiarity with [K97].

First of all, recall how the Kontsevich deformation quantization
formula is written. Let $\alpha$ be a Poisson structure on $V^*$.
Then the formula is
\begin{equation}\label{eq2.7.1}
f\star g=f\cdot g+\sum_{k\ge 1}\hbar^k(\sum_{m\ge 1}\frac
1{m!}\sum_{\Gamma \in G_{2,m}^2}W_\Gamma
U_\Gamma(\alpha,\dots,\alpha))
\end{equation}
Here $\Gamma$ is an admissible graph with two vertices on the "real
line" and $m$ vertices in the upper half-plane, and with 2 outtgoing
edges at each vertex in the upper half-plane, that is, $\Gamma\in
G_{2,m}^2$, in particular, it is an oriented graph with $2m$ edges;
$W_\Gamma$ is the Kontsevich integral of the Graph $\Gamma$. {\it
Let us note that the all graphs involved in (\ref{eq2.7.1}) may have
arbitrary many incoming edges at each vertex at the upper
half-plane, and exactly two outgoing edges}.

Now let $\alpha=\sum_{ij}\alpha_{ij}\partial_i\wedge \partial_j$,
where $\alpha_{ij}=\sum_I c_{ij}^I x_{i_1}\dots x_{i_k}$ ($I$ is a
multi-index).

Then the "Koszul dual" polyvector $\K(\alpha)$ is a polyvector field
with quadratic coefficients:

\begin{equation}\label{eq2.7.2}
\K(\alpha)=\sum_{i,j,I}c_{ij}^I (\xi_i\xi_j) \cdot
\partial_{\xi_{i_1}}\wedge\dots\wedge\partial_{\xi_{i_k}}
\end{equation}
It has total degree 1 and satisfies the Maurer-Cartan equation.

Firstly we write the formula for the image of $\K(\hbar\alpha)$ by
the Kontsevich formality, that is, denote by
\begin{equation}\label{eq2.7.3}
\mathcal{U}(\K(\alpha))=\hbar U_1(\K(\alpha))+\hbar^2
\frac12U_2(\K(\alpha),\K(\alpha))+\dots+\hbar^k\frac1{k!}U_k(\K(\alpha),\dots,\K(\alpha))+\dots
\end{equation}
We can write down explicitly this formula in graphs. {\it Let us
note the the graphs involving in (\ref{eq2.7.3}) may have arbitrary
many outgoing edges in the vertices at the upper half-plane, but
exactly two, one, or 0 incoming edges, because all components of
$\K(\alpha)$ are quadratic polyvector fiels}. That is, in a sense
the graphs in (\ref{eq2.7.3}) are "dual" to the graphs in
(\ref{eq2.7.1}).

Let us note also that the right-hand side of (\ref{eq2.7.3}) is a
polydifferential operator in $\Hoch(\Lambda(V^*))$ of
non-homogeneous Hochschild degree, but of the total (Hochschild
degree and $\Lambda$-degree) +1.

Now we should apply to $\mathcal{U}(\K)$ our map $\Phi_1$ to get a
derivation of the cobar-complex $CoBar(\Lambda^-(V))$. After this,
we get the final answer for $d_\hbar$.

Let us compute its component in the first power of $\hbar$. It is
just the Hochschild-Kostant-Rosenberg map of $U_1(\K(\alpha))$ which
is the symmetrization map $\Sym\colon S(V)\to T(V)$ in this case.

Let us note that in the case of a quadratic Poisson structure
algebra the relations $R_{ij}\in (V\otimes V)\xxx\subset T(V)\xxx$
are quadratic.

One can imagine then the relations $x_i\otimes x_j-x_j\otimes
x_i=d_\hbar(\xi_i\wedge \xi_j)$ hold {\it exactly} in the Kontsevich
star-algebra defined from the same propagator. It would be a "pure
duality", which would be very nice. Our Conjecture 2.6 is a weaker
statement, that this relations hold in an algebra gauge equivalent
to the Kontsevich star-algebra. Anyway, these relations may
considered as some very non-trivial relations with Kontsevich
integrals, which are transcendental numbers which in general is
almost impossible to compute directly.

\subsection*{Acknowledgements}
I am thankful to Misha Bershtein, Pavel Etingof, Borya Feigin,
Giovanni Felder, and to Bernhard Keller for interesting and useful
discussions. The work was partially supported by the research grant
R1F105L15 of the University of Luxembourg.

Faculty of Science, Technology and Communication, Campus
Limpertsberg, University of Luxembourg,
162A avenue de la Faiencerie, L-1511 LUXEMBOURG\\
{\it e-mail}: {\tt borya$\_$port@yahoo.com}

\end{document}